\documentclass{article}

\usepackage{arxiv}

\usepackage[utf8]{inputenc} 
\usepackage[T1]{fontenc}    
\usepackage{hyperref}       
\usepackage{url}            
\usepackage{booktabs}       
\usepackage{amsfonts}       
\usepackage{nicefrac}       
\usepackage{microtype}      
\usepackage{lipsum}		
\usepackage{graphicx}
\usepackage{natbib}
\usepackage{doi}

\usepackage{soul}
\usepackage[utf8]{inputenc}
\usepackage{graphicx}
\usepackage{amsmath}
\usepackage{amsthm}
\usepackage{booktabs}
\usepackage{txfonts}
\usepackage{epsfig}

\usepackage{amssymb}
\usepackage{epsfig}
\usepackage{mathtools}
\usepackage[boxruled]{algorithm2e}

\DeclareMathOperator*{\argmin}{\arg\!\min}
\usepackage{multirow}
\usepackage{booktabs}
\usepackage{balance}
\usepackage{subfig}

\newtheorem{theorem}{Theorem}
\newtheorem{definition}{Definition}
\newtheorem{assumption}{Assumption}

\DeclareMathOperator{\Tr}{Tr}

\newcommand{\norm}[1]{\left\lVert#1\right\rVert}

\usepackage[colorinlistoftodos]{todonotes}

\title{Causal Inference in Network Economics}

\date{} 					

\author{

       Sridhar Mahadevan\\
   Adobe Research, 345 Park Avenue, San Jose, CA 95110 \\
   
   smahadev@adobe.com
}

\begin{document}
\maketitle

\begin{abstract}
Network economics is the study of a rich class of equilibrium problems that occur in the real world, from traffic management to supply chains and two-sided online marketplaces. 
 In this paper we explore causal inference in network economics, building on the mathematical framework of  {\em variational inequalities}, which is a generalization of classical optimization. Our framework can be viewed as a synthesis of the well-known variational inequality formalism with the broad principles of  causal inference. 
\end{abstract}

\section{Introduction}\label{sec:intro}

\begin{figure}[h]
\begin{center}
\begin{minipage}[t]{0.45\textwidth}
\includegraphics[width=\textwidth,height=1.75in]{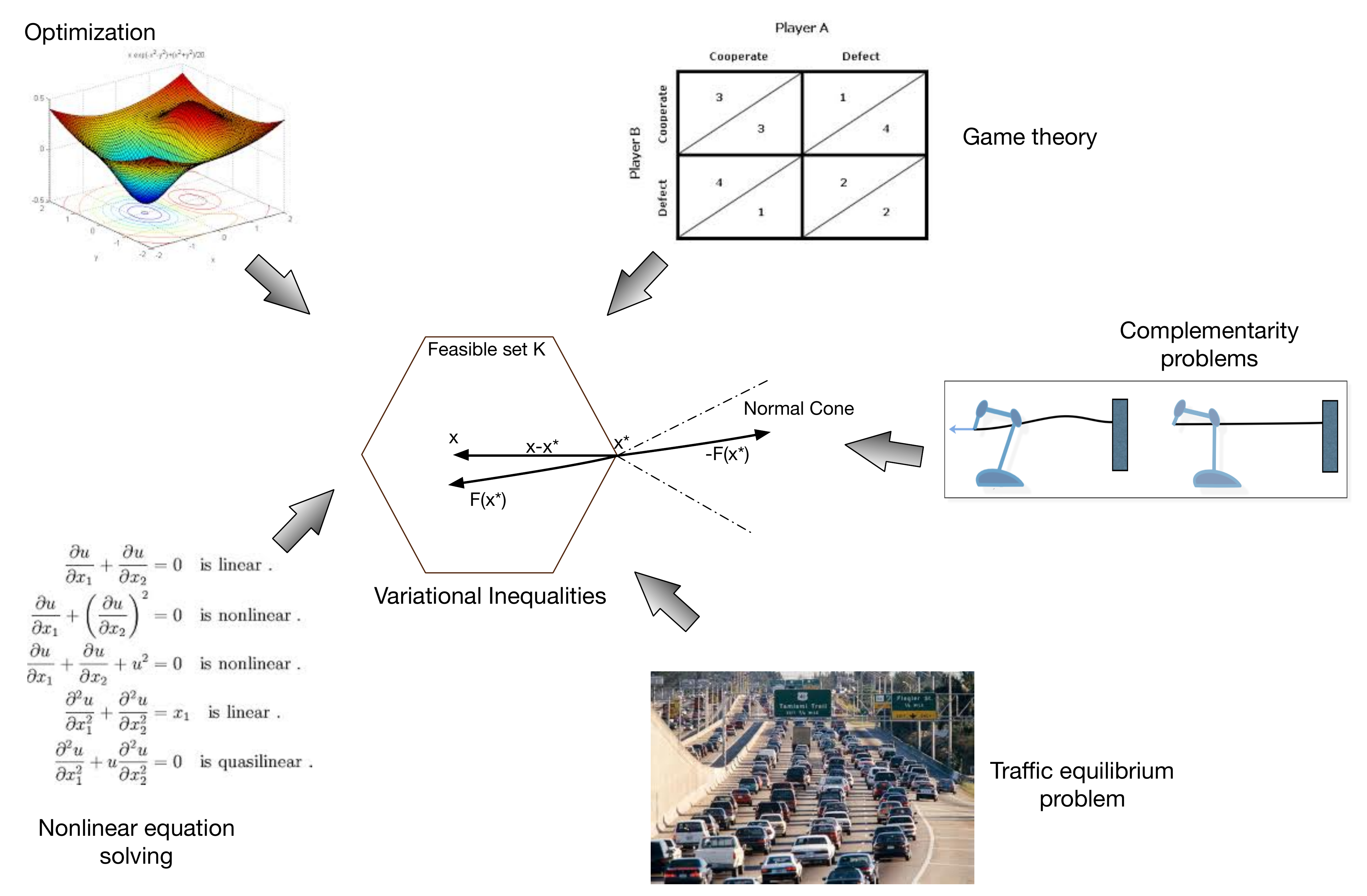}
\end{minipage}
\end{center}
\caption{Applications of causal variational inequalities.}
\label{vi-applications}
\end{figure}

We explore a variational calculus for causal inference, specifically our approach combines two successful, but hitherto largely distinct, areas: variational inequalities \citep{facchinei-pang:vi,stampacchia}, a mathematical formalism developed in physics in the mid-1960s for modeling equilibrium problems in continuum mechanics, but which generalizes (convex) optimization, non-cooperative games, fixed point equation solving, and complementarity problems; and the broad family of  causal inference models, which includes path diagrams, potential outcomes, structural equations, and graphical models \citep{pearl-book,rubin-book}. VI's have been applied to a wide range of problems in network economics \citep{nagurney:vibook}, game theoretic generative adversarial networks \citep{gemp-gan}, and reinforcement learning \citep{bertsekas:vi-tr09,rl-lcp} (see Figure~\ref{vi-applications}); 

Our proposed framework of causal variational inequalities extends the classical variational inequality model, which is usually defined as a set of deterministic multidimensional vector field mappings $F_i: K_i \rightarrow \mathbb{R}^n$, where $K_i$ is a convex set defined over a set of exogenous non-manipulable variables $U$ and a set of endogenous modifiable variables $V$. The ensemble of vector field mappings jointly defines a global vector field mapping $F: K \rightarrow \mathbb{R}^n$. Solving a variational inequality means finding an element $x \in K$ where the set of inequalities holds: 
\begin{equation}
 \langle F(x^*), (x - x^*) \rangle \geq 0, \ \ \forall x \in K
 \end{equation} 
where $\langle .,.\rangle$ denotes the inner product in $\mathbb{R}^n$. \footnote{Variational inequalities were originally developed for infinite-dimensional Hilbert spaces.  We consider only $n$-dimensional Euclidean spaces.} To see the connection with optimization, note that if $F(x) = \nabla f(x)$, the gradient of a (convex) multi-dimensional function, then the above condition is precisely the requirement for $x^*$ to be the (global) minimum. Broadly speaking, we will interpret causal intervention as inducing a submodel mapping $F_w$, which usually may be due to the manipulation of some endogenous variable $w$. However, our framework is agnostic on the particular manipulation mechanisms, and extends to other possibilities, such as modifying a pricing function or a change in demand. As a simple example to build intuition, consider a network economy consisting of $m$ producers of personal protective equipment (PPE), and $n$ demand markets. For trade to occur between the producer $i$ and the demand market $j$,  the following equilibrium condition must be satisfied: 
\begin{eqnarray}
\label{eqcond}
\pi_i(Q) + c_{ij}(Q) &=& \rho_j \ \mbox{if} \ \ Q^{*}_{ij} > 0 \\
\pi_i(Q) + c_{ij}(Q) &\geq& \rho_j \ \mbox{if} \ \ Q^*_{ij} = 0
\end{eqnarray}
This condition asserts that  trade, measured by $Q_{ij}$,  between producer $i$ and demand market $j$ will occur precisely when the supply price $\pi_i$ plus the transportation cost $c_{ij}$ is equal to the demand price $\rho_j$, else $Q^*_{ij} = 0$. Note that supply price and transportation costs are a function of trade volume $Q$ over the whole network. We can write this equilibrium condition in the form of a variational inequality as follows: 
\[ (\pi_i(Q) + c_{ij}(Q) - \rho_j) (Q_{ij} - Q^*_{ij}) \geq 0 \]
Note that if $Q^*_{ij} > 0$, then according to Equation~\ref{eqcond}, the first term above must equal $0$. Alternatively, if $Q^*_{ij} = 0$, then the reason is that the combined cost of manufacturing $\pi_i$ and transportation $c_{ij}$ exceeds the demand price $\rho_j$, hence the inequality is again satisfied. We can view each equilibrium condition above as a submodel $F_{ij}(Q)$, and there will be one such submodel for each producer and consumer. The number of submodels is equal to the number of trade paths, which collectively define the overall vector field mapping $F$. We can now model causal interventions in this system, such as raising or lowering prices, as inducing a modified model $F_w$, and analyze the effect of these causal interventions on the equilibrium trade between producers and suppliers. 
 
In the remainder of the paper, we define causal variational inequalities more precisely, and characterize its theoretical properties. We illustrate our variational causal calculus with a range of examples, including a multi-tiered complex network economics non-cooperative game involving a group of ``producer" agents whose goal is to sell content (digital media, physical goods, or financial instruments) to a set of demand markets (users), where the transportation of the goods is under the control of transport network agents. Our proposed causal formulation enables us to  formulate  and study causal interventions in not only complex network games, but also  to many other applications of VI's. 

\section{Causal Variational Inequalities}
\label{cvi} 

Our variational formulation of causal inference  is a synthesis of classical variational inequalities \citep{facchinei-pang:vi} and causal models \citep{rubin-book,pearl-book}.  More precisely, a causal variational inequality model $\cal{M} =$ CVI($F,K$), where $F$ is a collection of modular vector-valued functions defined as $F_i$, where $F_i: K_i \subset \mathbb{R}^{n_i} \rightarrow \mathbb{R}^{n_i}$, with each $K_i$ being a convex domain such that $\prod_i K_i = K$. We assume that the domains of each $F_i$ range over a collection $V$ of endogenous variables, and a set $U$ of exogenous variables, where only the endogenous variables are subject to causal manipulation. We model each intervention as a submodel $F_w$, and each component of $F_w$ reflects the effect of some manipulation of a subset $V_w \subseteq V$ of endogenous variables. 

\begin{definition}
\label{cvidef}
The finite-dimensional causal  variational inequality problem is defined by a model ${\cal M}$ = CVI($F, K)$, where the vector-valued mapping $F$ depend on both deterministic and stochastic elements, namely $F(x) = E[F(x, \eta)]$. where $\eta$ is a random variable defined over the probability space $(\Omega, {\cal F}, P)$, $E[.]$ denotes expectation with respect to the probability distribution $P$ over the random variable $\eta$, and $F: K \rightarrow \mathbb{R}^n$ is a given continuous function, $K$ is a given closed convex set, and $\langle .,.\rangle$ is the standard inner product in $\mathbb{R}^n$. A causal intervention is modeled as a submodel ${\cal M}_w$ = CVI($F_w, K)$, where $F_w(x) = E_w[F(x, \eta | \hat{w})]$, where $\hat{w}$ denotes the intervention of setting of variable $w$ to a specific non-random value, and where $E_w[.]$ now denotes expectation with respect to the intervention probability distribution $P_w$. Solving a causal VI is defined as finding a vector $x^* = (x^*_1, \ldots, x^*_n) \in K \subset \mathbb{R}^n$ such that
\begin{equation*}
\langle F_w(x^*), (y - x^*) \rangle \geq 0, \ \forall y \in K
\end{equation*}
\end{definition}
\begin{figure}[h]
\begin{center}
\begin{minipage}[t]{0.45\textwidth}
\includegraphics[width=\textwidth,height=1.25in]{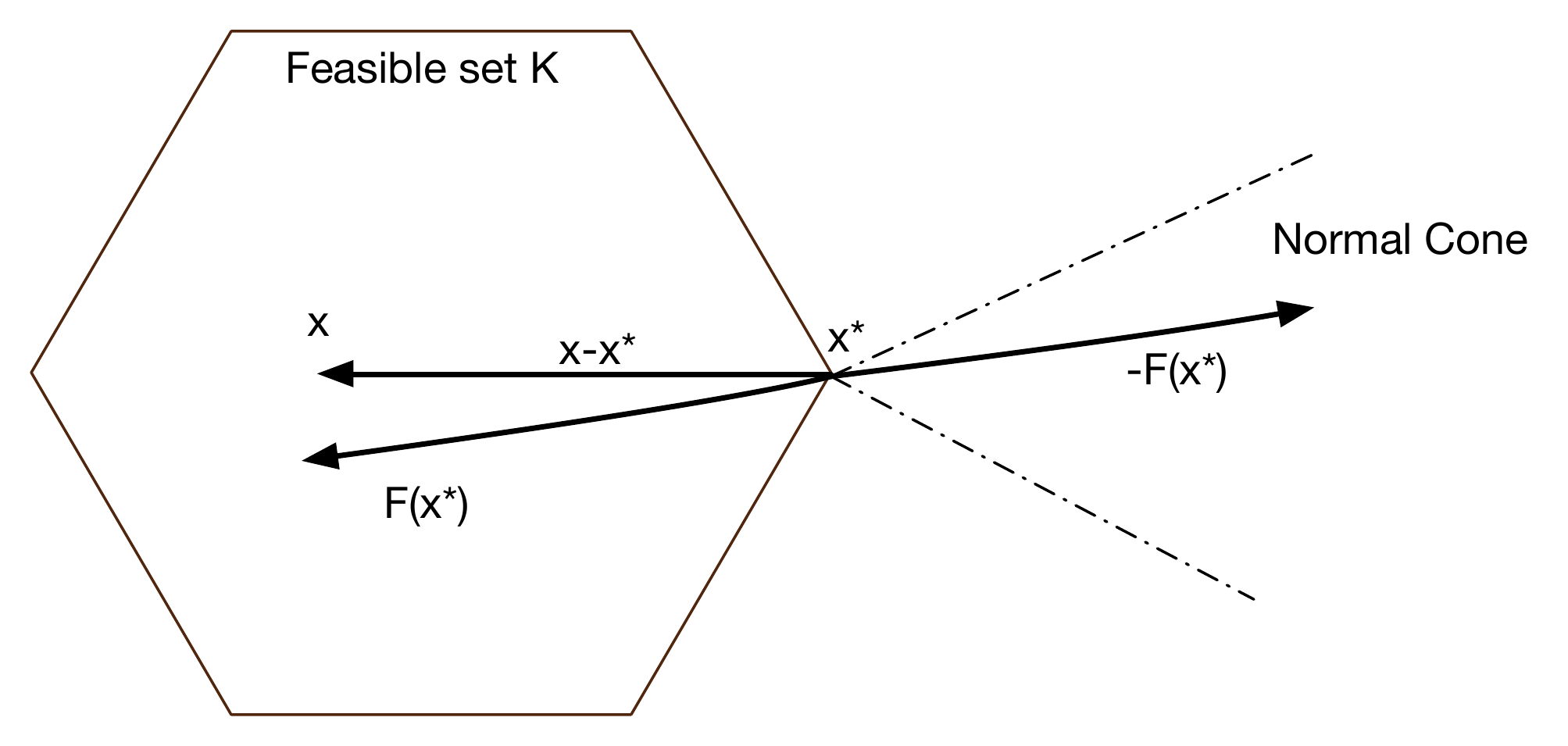}
\end{minipage}
\end{center}
\caption{This figure provides a geometric interpretation of a causal variational inequality $CVI(F_w,K)$. The mapping $F_w$ defines a vector field over the feasible set $K$ and a probability space, where $E_w(F(x, \eta) | \hat{w})$ is the conditional mean vector field (denoted in the figure by $F$), computed over the intervention distribution $P_w$. At the solution point $x^*$, the vector field $F(x^*)$ is directed inwards at the boundary, and  $-F(x^*)$ is an element of the normal cone $C(x^*)$ of $K$ at $x^*$ where the normal cone  $C(x^*)$ at the vector $x^*$ of a convex set $K$ is defined as $C(x^*) = \{y \in \mathbb{R}^n | \langle y, x - x^* \rangle \leq 0, \forall x \in K \}$.}
\label{vi-geom}
\end{figure}

\subsection{Causal Traffic Model}

\begin{figure}[h]
\begin{center}
\begin{minipage}[t]{0.35\textwidth}
\includegraphics[width=\textwidth,height=1.25in]{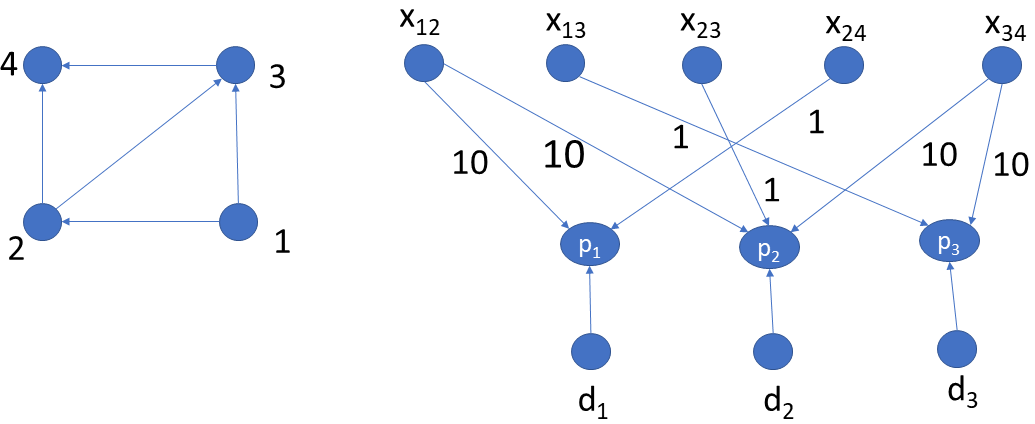}
\end{minipage}
\end{center}
\caption{A simple road network (left) and part of its associated structural causal model on the right.}
\label{traffic}
\end{figure}

Figure~\ref{traffic} illustrates the problem of determining the equilibrium flow on a small network of roads from source location $1$ to a destination location $4$. This example is an instance of the famous Braess paradox \citep{braess}.  Let us order the directed edges as $E = \{ (1,2), (1,3), (2,3), (2,4), (3,4) \}$. We are given that $6$ cars per unit time leave node $1$ intending to reach node $4$. How does the traffic distribute across the edges?  The vector field $F$ for this problem representing delays across the edges is as follows (with $d_1= d_3 = 50$, and $d_1=10$):\footnote{For simplicity, $F$ here is symmetric. Section~\ref{soi} will illustrate a non-symmetric $F$ mapping.} 
\begin{equation}
\label{braess}
F(x) = \left[ \begin{array}{ccccc}  10 & 0 & 0 & 0 & 0 \\ 0 & 1 & 0 & 0 & 0 \\ 0 & 0 & 1 & 0 & 0 \\ 0 & 0 & 0 & 1 & 0 \\ 0 & 0 & 0 & 0 & 10   \end{array} \right] \left[ \begin{array}{c} x_{12}  \\ x_{13} \\ x_{23} \\ x_{24} \\ x_{34} \end{array} \right] + \left[ \begin{array}{c} 0 \\ 50 \\ 10 \\ 50 \\ 0 \end{array} \right].
\end{equation}
The space of feasible solutions $K = \{ x \in \mathbb{R}_+^5 | \ B x = b \}$, where the $B$ is the incidence matrix for the graph, and the vector $b$ is given as: 
\begin{equation}
B = \left[ \begin{array}{ccccc} 1 & 1 & 0 & 0 & 0 \\ -1 & 0 & 1 & 1 & 0 \\ 0 & -1 & -1 & 0 & 1 \\ 0 & 0 & 0 & -1 & -1 \end{array} \right], \ b = \left[ \begin{array}{c} 6 \\ 0 \\ 0 \\ -6 \end{array} \right]. 
\end{equation}
Consider the causal intervention ${\bf do}(x_{23} = 0)$, in which case no cars flow through the diagonal edge from node 2 to 3. There are only two paths in this case: $p_1$ with edges (1,2), and (2,4), and $p_3$ with edges (1,3) and (3,4). The drivers will distribute along these two paths such that the travel times are equal, which implies $10 x_{12} + x_{24} + 50 = x_{13} + 50 + 10 x_{34}$. The equilibrium solution from solving the causal variational inequality $CVI(F_{\hat{x}_{23}},K)$ is given by $x^*_{\hat{x}_{23}} = (3,3,0,3,3)$ (flows across edges in $E$), giving a total delay of $83$ time units across either path. Consider now the original network as shown in Figure~\ref{traffic}, where traffic is allowed to flow through edge (2,3). Surprisingly, allowing traffic to flow across this diagonal edge increases delay time to 92, because the solution to the ``untreated" CVI(F,K) is $x^* = (4,2,2,2,4)$. It turns out that the Braess paradox is only realized for demands in a narrow range, and if demands fall below and above this range, the paradox disappears. For simplicity, we have only considered linear cost functions for delay, but our causal VI formulation extends to the case when the delay functions are nonlinear \citep{braess-2010}.
\subsection{Properties of Mappings} 
The solution to a (causal) VI depends on the properties satisfied by the mapping $F$ and the feasible space $K$. If $K$ is compact and $F$ is continuous, it is straightforward to prove using Brower's fixed point theorem that there is always at least one solution to any VI (see Theorem~\ref{projthm}).  However, to obtain a unique solution, a stricter condition is necessary. 
\begin{definition}
$F(x)$ is {\em monotone} if $\langle F(x) - F(y), x - y \rangle \ge 0$, $\forall x, y \in K$.
\end{definition}

\begin{definition}
$F(x)$ is {\em strongly monotone} if $\langle F(x) - F(y), x  - y \rangle \geq \mu \| x - y \|^2_2, \mu > 0, \forall x,y \in K$.
\end{definition}
\begin{definition}
$F(x)$ is {\em Lipschitz} if $\| F(x) - F(y) \|_2 \leq L \|x - y \|_2, \forall x,y \in K$.
\end{definition}

Crucially, VI problems can only be converted into equivalent optimization problems when a very restrictive condition is met on the Jacobian of the mapping $F$, namely that it be symmetric. Most often, real-world applications of VIs, such as the example in Section~\ref{soi}, do not induce symmetric Jacobians. 
\begin{theorem}
\label{equivalence}
Assume $F(x)$ is continuously differentiable on $K$ and that the Jacobian matrix $\nabla F(x)$ of partial derivatives of $F_i(x)$ with respect to (w.r.t) each $x_j$ is symmetric and positive semidefinite. Then there exists a real-valued convex function $f: K \rightarrow \mathbb{R}$ satisfying $\nabla f(x) = F(x)$ with $x^*$, the solution of  VI(F,K), also being the mathematical programming problem of minimizing $f(x)$ subject to $x \in K$.
\end{theorem}
The algorithmic development of methods for solving VIs begins with noticing their connection to fixed point problems.
\begin{theorem}
\label{projthm}
The vector $x^*$ is the solution of VI(F,K) if and only if, for any $\alpha > 0$, $x^*$ is also a fixed point of the map  $x^* = P_K(x^* - \alpha F(x^*))$,
where $P_K$ is the projector onto convex set $K$.
\end{theorem}
In terms of the geometric picture of a VI illustrated in Figure \ref{vi-geom}, this property means that the solution of a VI occurs at a vector $x^*$ where the vector field $F(x^*)$ induced by $F$ on $K$ is normal to the boundary of $K$ and directed inwards, so that the projection of $x^* - \alpha F(x^*)$ is the vector $x^*$ itself. This property forms the basis for the projection class of methods that solve for the fixed point.

\subsection{Projected Dynamical Systems}
\label{theory_pds}
Causal Projected dynamical systems (CPDS) provide an alternate perspective of VIs that's helpful in terms of analyzing stability of the equilibria and designing algorithms, and particularly useful for understanding the short-term and long-term effect of a causal intervention \citep{DBLP:conf/nips/ToulisP16}. Any equilibrium points of a CPDS(F,K) coincide with the solutions of a causal VI(F,K).  The principal difference between this formulation and that of ODEs is that the right hand side of the differential equation is discontinuous, to incorporate constraints such as non-negativity on the variables. The operator $\Pi_K(x,v)$ below denotes the Gateaux differential of vector $x$ along the direction $v$, and $P_K$ denotes the usual projection onto a convex space $K$. In a causal PDS, we want to model the change in solution resulting from causal interventions that define a modified vector field $F_w$: 
\begin{definition}
\label{pdsdef}
A causal projected dynamical system, CPDS($F_w,K$), corresponding to a causal variational inequality CVI$(F_w,K)$is defined as the set of differential equations: 
\begin{center}$\dot{X} = \Pi_{K}(X,-F_w(X))$  with  $X(0) = X^0$  and\end{center}
\[ \Pi_{K}(X,-F_w(X)) = \lim_{\delta \rightarrow 0} \frac{P_{K}(X-\delta F_w(X))-X}{\delta} \]
\end{definition}

\subsection{Variational Inequalities and Games}
\label{theory_nash}
Together, the causal VI and PDS frameworks provide a mathematically elegant approach to modeling and solving equilibrium problems in game theory \citep{fudenbergtheory,nisan2007algorithmic}.  A {\em Nash game} consists of $m$ players, where player $i$ chooses a strategy $x_i$ belonging to a closed convex set $X_i \subset \mathbb{R}^n$.  After executing the joint action, each player is penalized (or rewarded) by the amount $f_i(x_1,\ldots,x_m)$, where $f_i: \mathbb{R}^n \rightarrow \mathbb{R}$ is a continuously differentiable function.  A set of strategies $x^* = (x_1^*,\ldots,x_m^*) \in \Pi_{i=1}^M X_i$ is said to be in equilibrium if no player can reduce the incurred penalty (or increase the incurred reward) by unilaterally deviating from the chosen strategy.  If each $f_i$ is convex on the set $X_i$, then the set of strategies $x^*$ is in equilibrium if and only if $\langle \nabla_i f_i (x_i^*), (x_i - x_i^*) \rangle \ge 0$.  In other words, $x^*$ needs to be a solution of the VI $\langle F(x^*), (x-x^*) \rangle \ge 0$, where $F(x) = ( \nabla f_1(x), \ldots, \nabla f_m(x))$.

Two-person Nash games are closely related to {\em saddle point} optimization problems \citep{juditsky2011first,juditsky2011second,liu2012regularized} where we are given a function $f: X \times Y \rightarrow \mathbb{R}$, and the objective is to find a solution $(x^*,y^*) \in X \times Y$ such that \begin{equation} f(x^*,y) \le f(x^*,y^*) \le f(x,y^*), \forall x \in X, \forall y \in Y. \end{equation}  Here, $f$ is convex in $x$ for each fixed $y$, and concave in $y$ for each fixed $x$. 

Complementarity problems provide the foundation for a number of Nash equilibrium algorithms.  The class of complementarity problems can also be reduced to solving a VI. When the feasible set $K$ is a cone, meaning that if $x \in K$, then $\alpha x \in K, \alpha \geq 0$, then the VI becomes a CP.
\begin{definition}
Given a cone $K \subset \mathbb{R}^n$ and mapping $F: K \rightarrow \mathbb{R}^n$, the complementarity problem CP(F,K) is to find an $x \in K$ such that $F(x) \in K^*$, the dual cone to $K$, and $\langle F(x), x \rangle \geq 0$. \footnote{Given a cone $K$, the dual cone $K^*$ is defined as $K^* = \{ y \in \mathbb{R}^n | \langle y, x \rangle \geq 0, \forall x \in K \}$.}
\end{definition}

The nonlinear complementarity problem (NCP) is to find $x^* \in \mathbb{R}^n_+$ (the non-negative orthant) such that $F(x^*) \geq 0$ and $\langle F(x^*), x^* \rangle = 0$. The solution to an NCP and the corresponding $VI(F, \mathbb{R}^n_+)$ are the same, showing that NCPs reduce to VIs. In an NCP, whenever the mapping function $F$ is affine, that is $F(x) = Mx + b$, where $M$ is an $n \times n$ matrix, the corresponding NCP is called a linear complementarity problem (LCP) \citep{murty:lcpbook}.

\section{Causal Network Economics}
\label{soi}

\begin{figure}[h!]
\centering
\begin{minipage}{1\textwidth}
\includegraphics[scale=0.6]{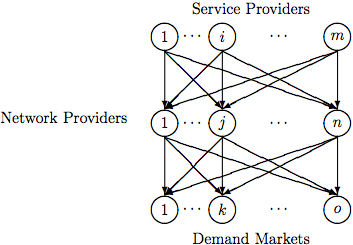}
\end{minipage}
\caption{A generic network economic model based on \citep{nagurney:soi}.  Demand markets consisting of individual users or groups of users, who choose  combinations of service providers and transport providers, both of whom compete with each other for demand markets.}
\label{SOI-diagram}
\end{figure}

We now describe how to model causal inference in a network economics problem, which will be useful in illustrating the abstract definitions from the previous section. Network economics games can be considerably more complex than single-layer games used in some previous work on causal inference  \citep{DBLP:journals/jmlr/BottouPCCCPRSS13,DBLP:conf/nips/ToulisP16,Wager2019ExperimentingIE}. The model in Figure~\ref{SOI-diagram} is drawn from \citep{nagurney:vibook,nagurney:soi}. which were  deterministic, and included no analysis of causal interventions. This network economics model comprises of three tiers of agents: producer agents, who want to sell their goods, transport agents who ship merchandise from producers, and demand market agents interested in purchasing the products or services. The model applies both to electronic goods, such as video streaming, as well as physical goods, such as face masks and other PPEs. 

The model assumes $m$ service providers, $n$ network providers, and $o$ demand markets.  Each firm's utility function is defined in terms of the nonnegative service quantity (Q), quality (q), and price ($\pi$) delivered from service provider $i$ by network provider $j$ to consumer $k$.  Production costs, demand functions, delivery costs, and delivery opportunity costs are designated by $f$, $\rho$, $c$, and $oc$ respectively.  Service provider $i$ attempts to maximize its utility function $U_i^1(Q,q^*,\pi^*)$ by adjusting $Q_{ijk}$ (eqn. \ref{U1}).  Likewise, network provider $j$ attempts to maximize its utility function $U_j^2(Q^*,q,\pi)$ by adjusting $q_{ijk}$ and $\pi_{ijk}$ (eqn. \ref{U2}).

\begin{subequations}
\begin{align}
\label{U1}
U_i^1(Q,q^*,\pi^*) &= \sum_{j=1}^n \sum_{k=1}^o \hat{\rho}_{ijk}(Q,q^*)Q_{ijk} - \hat{f}_i(Q)\\
&- \sum_{j=1}^n \sum_{k=1}^o \pi^*_{ijk}Q_{ijk}, \hspace{0.2cm} Q_{ijk} \ge 0 \nonumber
\end{align}
\begin{align}
\label{U2}
U_j^2(Q^*,q,\pi) = &\sum_{i=1}^m \sum_{k=1}^o \pi_{ijk}Q^*_{ijk}\\
- &\sum_{i=1}^m \sum_{k=1}^o (c_{ijk}(Q^*,q) + oc_{ijk}(\pi_{ijk})), \nonumber \\
&q_{ijk}, \pi_{ijk} \ge 0 \nonumber
\end{align}
\end{subequations}

We assume the governing equilibrium is Cournot-Bertrand-Nash and the utility functions are all concave and fully differentiable.  This establishes the equivalence between the equilibrium state we are searching for and the variational inequality to be solved where the $F$ mapping is a vector consisting of the negative gradients of the utility functions for each firm.  Since $F$ is essentially a concatenation of gradients arising from multiple independent, conflicting objective functions, it does not correspond to the gradient of any single objective function.  

\begin{subequations}
\begin{align}
\label{SOI-vi}
&\langle F(X^*),X-X^* \rangle \ge 0, \forall X \in \mathcal{K},\\
\text{where } \hspace{0.2cm} &X = (Q,q,\pi) \in \mathbb{R}^{3mno+} \nonumber \\
\text{and } \hspace{0.35cm} &F(X) = (F^1_{ijk}(X), F^2_{ijk}(X), F^3_{ijk}(X)) \nonumber
\end{align}
\begin{align}
F^1_{ijk}(X) &= \frac{\partial f_i (Q)}{\partial Q_{ijk}} + \pi_{ijk} - \rho_{ijk} - \sum_{h=1}^n \sum_{l=1}^o \frac{\partial \rho_{ihl} (Q,q) }{\partial Q_{ijk}} \times Q_{ihl} \label{F1} \\
F^2_{ijk}(X) &= \sum_{h=1}^m \sum_{l=1}^o \frac{\partial c_{hjl} (Q,q)} {\partial q_{ijk}} \label{F2} \\
F^3_{ijk}(X) &= -Q_{ijk} + \frac{\partial oc_{ijk}(\pi_{ijk})}{\partial \pi_{ijk}} \label{F3}
\end{align}
\end{subequations}

The variational inequality in Equations~\ref{SOI-vi} represents the result of combining the utility functions of each firm into standard form.  $F^1_{ijk}$ is derived by taking the negative gradient of $U_i^1$ with respect to $Q_{ijk}$.  $F^2_{ijk}$ is derived by taking the negative gradient of $U_j^2$ with respect to $q_{ijk}$.  And $F^3_{ijk}$ is derived by taking the negative gradient of $U_j^2$ with respect to $\pi_{ijk}$.

\subsection{Numerical Example} 

We extend the simplified numerical example in \citep{nagurney:soi} by adding stochasticity to illustrate our causal variational formalism. Let us assume that there are two service providers, one transport agent, and two demand markets. Define the production cost functions:
\[ f_1(Q) = q^2_{111} + Q_{111}  + \eta_{f_1}, f_2(Q) = 2 Q^2_{111} + Q_{211} + \eta_{f_2} \] 
where $\eta_{f_1}, \eta_{f_2}$ are random variables indicating errors in the model. 
Similarly, define the demand price functions as:
\begin{eqnarray*}
 \rho_{111}(Q,q) = -Q_{111} - 0.5 Q_{211} + 0.5 q_{111} + 100 + \eta_{\rho_{111}} \\
 \rho_{211}(Q,q) = -Q_{211} - 0.5 Q_{111} + 0.5 q_{211} + 200 + \eta_{\rho_{211}}
 \end{eqnarray*}
Finally, define the transportation cost functions as:
\begin{eqnarray*}
c_{111}(Q,q) = 0.5(q_{111} - 20)^2 + \eta_{c_{111}} \\ 
c_{211}(Q,q) = 0.5(q_{211} - 10)^2 + \eta_{c_{211}} 
\end{eqnarray*}
and the opportunity cost functions as:
\[ oc_{111}(\pi_{111}) = \pi_{111}^2 + \eta_{oc_{111}}, oc_{211}(\pi_{211}) =\pi_{211}^2 + \eta_{oc_{211}} \] 
Using the above equations, we can easily compute the component mappings $F_i$ as follows: 
\begin{eqnarray*}
F^1_{111}(X) = 4 Q_{111} + 0.5 Q_{211} - 0.5 q_{111} - 99 \\
F^1_{211}(X) = 6 Q_{211} + \pi_{211} - 0.5 Q_{111} - 0.5 q_{211} -199\\
F^2_{111}(X) = q_{111} - 20, \ F^2_{211}(X) = q_{211} - 10 \\
F^3_{111}(X) = -Q_{111} + 2 \pi_{111}, \ F^3_{211}(X) = -Q_{211} + 2 \pi_{211}
\end{eqnarray*}
For simplicity, we have not indicated the noise terms above, but assume each component mapping $F_i$ has an extra noise term $\eta_i$. It is also clear that we can now give precise semantics to causal intervention in this system, following the principles laid out in \citep{pearl-book}. For example, if we set the network service cost $q_{111}$ of network provider 1 serving the content producer 1 to destination market 1 to 0, then the production cost function under the intervention distribution is given by 
\[ E_{q_{111} = 0}(f_1(Q) | \hat{q}_{111}) = Q_{111} + E_{q_{111}=0} (\eta_{f_1} | \hat{q}_{111})\] 
Finally, the Jacobian matrix associated with $F(X)$ is given by the partial derivatives of each $F_i$ mapping with respect to $(Q_{111}, Q_{211}, q_{111}, q_{211}, \pi_{111}, \pi_{211})$ is given as:
\[ - \nabla U(Q,q,\pi) =  \left( \begin{array}{cccccc} 4 & .5 & -.5 & 0 & 1 & 0 \\
0.5 & 6 & 0 & -.5 & 0 & 1 \\ 0 & 0 & 1 & 0 & 0 & 0 \\ 0 & 0 & 0 & 1 & 0 & 0 \\ -1 & 0 & 0 & 0 & 2 & 0 \\ 0 & -1 & 0 & 0 & 0 & 2\end{array} \right) \]
Note this Jacobian is non-symmetric, but positive definite, as it is diagonally dominant. Hence the induced vector field $F$ can be shown to be strongly monotone, and the induced VI has exactly one solution. 

\section{Algorithms} 

 We now discuss algorithms for solving causal VI's. There are a wealth of existing methods for deterministic VI's \citep{facchinei-pang:vi,nagurney:vibook}), which can be adapted to solving causal VI's.   The simplest method for solving a causal VI is the well-known projection algorithm \citep{facchinei-pang:vi}: 
\[ x_{k+1} = \Pi_K [x_k - \alpha_k F_w(x_k)] \]
where $F_w$ is the vector field induced by some causal intervention, 
which can be viewed as a modification of the classical projection method for deterministic VI's. The algorithm follows the direction of the negative vector field at a point $x_k$, and if the iterate falls outside the feasible space $K$, it projects back into $K$. If $F_w$ is strongly monotone, and Lipschitz, and the learning rate $\alpha_k$ is suitably designed, then the projection algorithm is guaranteed to find the solution to a causal VI.

Understanding the convergence of the projection method will give us insight into how to analyze causal interventions in VI's. At the heart of convergence analysis of any VI method is bounding the iterates of the algorithm. In the below derivation, $x^*$ represents the final solution to a causal VI, and $x_{k+1}$, $x_k$ are successive iterates: 
\begin{small}
\begin{eqnarray*}
\| x_{k+1} - x^* \|^2 &=& \| P_K[x_k - \alpha_k F_w(x_k)] - P_K[x^* - \alpha_k F_w(x^*)] \|^2 \\
&\leq& \| (x_k - \alpha_k F_w(x_k)) - (x^* - \alpha_k F_w(x^*)) \|^2 \\
&=& \| (x_k - x^*) - \alpha_k (F_w(x_k) - F_w(x^*)) \|^2 \\
&=& \|x_k - x^* \|^2 - 2 \alpha_k \langle (F_w(x_k) - F_w(x^*)), x_k - x^* \rangle \\ &+& \alpha_k^2 \|F_w(x_k) - F_w(x^*)\|^2 \\
&\leq& (1 - 2 \mu \alpha_k + \alpha_k^2 L^2) \| x_k - x^*\|^2
\end{eqnarray*} 
\end{small}
Here, the first inequality follows from the nonexpansive property of projections, and the last inequality follows from strong monotonicity and Lipschitz property of the $F_w$ mapping. Bounding the term $\langle (F_w(x_k) - F_w(x^*)), x_k - x^* \rangle$ is central to the design of any VI method. As we show in the next section, in modeling causal interventions a similar term will arise, except under different mappings, representing the ``untreated" and "treated" cases. 

\citet{korpelevich} extended the projection algorithm with the well-known ``extragradient" method, which requires two projections, but is able to solve VI's for which the mapping $F$ is only monotone. \citet{GempM015} proposed a more sophisticated extragradient algorithm, combining Runge-Kutta methods from ODE's with a modified dual-space mirror-prox method \citep{nemirovski:siam,mirror-prox} to solve large network games modeled as VI's, such as the network economy described in Section~\ref{soi}, but required multiple projections corresponding to the order of the Runge-Kutta approximation. If projections are expensive, particularly in large network economy models, these algorithms may be less attractive than incremental stochastic projection methods, which we turn to next. 

\subsection{Incremental Projection Methods} 

We now describe an incremental two-step projection method for solving causal VI's, based on work by \citet{DBLP:journals/mp/WangB15}. Their algorithm adapted to causal VI's can be written as follows:
\begin{equation}
\label{2stepalg}
z_k = x_k - \alpha_k F_w(x_k, v_k), \ \ \ x_{k+1} = z_k - \beta_k (z_k - P_{w_k} z_k)
\end{equation}
where $\{v_k\}$ and $\{w_k\}$ are sequences of random variables, generated by sampling the causal VI model, and $\{\alpha_k\}$ and $\{\beta_k\}$ are sequences of positive scalar step sizes. Note that an interesting feature of this algorithm is that the sequence of iterates $x_k$ is not guaranteed to remain within the feasible space $K$ at each iterate. Indeed, $P_{w_k}$ represents the projection onto a randomly sampled constraint $w_k$. 

The analysis of convergence of this algorithm is somewhat intricate, and we will give the broad highlights as it applies to causal VI's. Define the set of random variables ${\cal F}_k = \{v_0, \ldots, v_{k-1}, w_0, \ldots, w_{k-1}, z_0, \ldots, z_{k-1}, x_0, \ldots, x_k \}$. Similar to the convergence of the projection method, it is possible to show that the error of each iteration is {\em stochastically contractive}, in the following sense: 
\[ E[ \| x_{k+1} - x^* \|^2 | {\cal F}_k ] \leq (1 - 2 \mu_k \alpha_k + \delta_k) \| x_k - x^* \|^2 + \epsilon_k, \ \ \ \mbox{w. p. 1} \]
where $\delta_k, \epsilon_k$ are positive errors such that $\sum_{k=0}^\infty \delta_k < \infty$ and $\sum_{k=0}^\infty \epsilon_k < \infty$. The convergence of this method rests on the following principal assumptions, stated below: 
\begin{assumption}
The intervened causal mapping $F_w$ is strongly monotone, and the sampled mapping $F_w(.,v)$ is {\em stochastically Lipschitz continuous} with constant $L > 0$, namely:
\begin{equation}
 E[ \| F_w(x, v_k) - F_w(y, v_k) \|^2 | {\cal F}_k] \leq L^2 \| x - y \|^2. \forall x, y \in \mathbb{R}^n
 \end{equation}
\end{assumption}

\begin{assumption}
The intervened mapping $F_w$ is bounded, with constant $B > 0$ such that 
\begin{equation}
\| F_w(x^*) \| \leq B, \ \ \ E[\| F_w(x^*, v) \|^2 | {\cal F}_k] \leq B^2, \ \ \forall k \geq 0
\end{equation}
\end{assumption}

\begin{assumption}
The distance between each iterate $x_k$ and the feasible space $K$ reduces ``on average", namely: 
\begin{equation}
\| x - P_K \|^2 \geq \eta \max_{i \in M} \|x - P_{K_i} x\|^2
\end{equation}
where $\eta > 0$ and $M = \{1, \ldots, m \}$ is a finite set of indices such that $\prod K_i = K$. 
\end{assumption}
\begin{assumption}
\[ \sum_{k=0}^\infty \alpha_k = \infty, \sum_{k=0}^\infty \alpha_k^2 < \infty, \sum_{k=0}^\infty \frac{\alpha_k^2}{\gamma_k} < \infty \]
where $\gamma_k = \beta_k (2 - \beta_k)$. 
\end{assumption}
\begin{assumption}
Supermartingale convergence theorem: Let ${\cal G}_k$ denote the collection of nonnegative random variables $\{y_k\}, \{u_k \}, \{a_k\}, \{b_k\}$ from $i=0, \ldots, k$
\begin{equation}
E[y_{k+1} | {\cal G}_k] \leq (1 + a_k) y_k - u_k + b_k, \ \ \forall k \geq 0, \ \mbox{w.p. 1}
\end{equation}
and $\sum_{k=0}^\infty a_k < \infty$ and $\sum_{k=0}^\infty b_k < \infty$ w.p. 1. Then, $y_k$ converges to a nonnegative random variable, and $\sum_{k=0}^\infty u_k < \infty$.
\end{assumption}
\begin{assumption}
The random variables $w_k, k=0, \ldots$ are such that for $\rho \in (0, 1]$
\[ \inf_{k\geq 0} P(w_k = X_i | {\cal F}_k) \geq \frac{\rho}{m}, \ \ i=1, \ldots, m, \ \mbox{w.p.1} \]
namely, the individual constraints will be sampled sufficiently. Also, the sampling of the stochastic components $v_k, k=0, \ldots$ ensures that
\[ E[F_w(x_k, v_k) | {\cal F}_k] = F_w(x), \ \ \forall x \in \mathbb{R}^n, \ k \geq 0 \]
\end{assumption}

Given the above assumptions, it can be shown that two-step stochastic algorithm given in Equation~\ref{2stepalg} converges to the solution of a causal VI, namely:
\begin{theorem}
Given a finite-dimensional causal  variational inequality problem is defined by a model ${\cal M}$ = CVI($F, K)$, and a causal intervention, defined by the submodel ${\cal M}_w$ = CVI($F_w, K)$, where $F_w(x) = E_w[F(x, \eta | \hat{w})]$, where $\hat{w}$ denotes the intervention of setting of variable $w$ to a specific non-random value, and where $E_w[.]$ now denotes expectation with respect to the intervention probability distribution $P_w$, the two-step algorithm given by Equation~\ref{2stepalg} produces a sequence of iterates $x_k$ that converges almost surely to $x^*$, where 
\begin{equation*}
\langle F_w(x^*), (y - x^*) \rangle \geq 0, \ \forall y \in K
\end{equation*}
\end{theorem}
{\bf Proof:} The proof of this theorem largely follows the derivation given in \citep{DBLP:journals/mp/WangB15}, where the only difference is that in a causal VI problem, we are conditioning the stochastic VI on the intervention distribution $P_w$. $\qed$

Crucially, the efficiency of this two-step method depends on how the constraints $w_k$ are sampled, i.e. either randomly or in some order. As we discuss next, since our goal is to understand the impact of causal interventions, we can use the property of locality of interventions to streamline this two-step procedure. 

\section{Analyzing  Interventions}  
\label{theory}

A fundamental question in causal inference is to measure the effect of some intervention, by comparing potential outcomes across the ``treated" units with the ``untreated" units \citep{rubin-book}. We characterize treatment effects in causal variational inequalities under interventions, building on the existing results on sensitivity analysis of classical variational inequalities \citep{nagurney:vibook}. First, we characterize a causal irrelevance theorem for causal variational inequalities. 

\begin{theorem}
If $Y$ is probabilistically causally irrelevant to $X$, given $Z$, then  CVI$(F_{\hat{y},\hat{z}}(x), K)$ has the same solution as CVI$(F_{\hat{y}',\hat{z}}(x), K)$. 
\end{theorem}

{\bf Proof:} The proof is straightforward given the axioms of causal irrelevance \citep{pearl-book}. If $Y$ is causally irrelevant to $X$ given $Z$, then it follows that $P(x | \hat{y}, \hat{z}) = P(x | \hat{y}', \hat{z})$ for all $y, y', x, z$, namely if $\hat{z}$ is fixed, then changing the value of $y$ has no influence on the distribution of $x$. In this case, the mapping $F$ under the two intervention distributions remains identical. $\qed$

Now we examine the case when interventions do alter the solution to a causal VI, where our goal is to measure the change in solution in terms of properties of the ``untreated mapping $F_0$ and the ``treated" mapping $F_1$. 
\begin{theorem}
\label{sens-thm1}
Let the solution of the original ``untreated" CVI$(F_0,K)$be denoted by $x_0$, where $F_0$ is assumed to be strongly monotone, and (stochastically) Lipschitz, with $\mu$ being the coefficient in the strong monotonicity property. Given a causal intervention, the ``treated" CVI$(F_1, K)$ results in the modified solution vector $x_1$. Then it follows that
\begin{equation}
    \| x_1 - x_0 \| \leq \frac{1}{\mu} \| F_1(x_1) - F_0(x_1) \|
\end{equation}
\end{theorem}
{\bf Proof:} Since $x_0$ and $x_1$ solve the ``untreated" and "treated" causal VI's, respectively, it must follow that:
\begin{eqnarray*}
\langle F_0(x_0), y  - x_0 \rangle \geq 0, \ \ \forall y \in K \\
\langle F_1(x_1), y - x_1 \rangle \geq 0, \ \ \forall y \in K
\end{eqnarray*}
Substituting $y = x_1$ in the first equation above, and $y = x_0$ in the second equation, it follows that:
\[ \langle (F_1(x_1) - F_0(x_0)) , x_1 - x_0 \rangle \leq 0 \] 
Equivalently, we get
\[ \langle (F_1(x_1) - F_0(x_0) + F_0(x_1) - F_0(x_1)), x_1 - x_0 \rangle \leq 0 \]
Using the monotonicity property of $F_0$, we get:
\begin{eqnarray*}
\langle (F_1(x_1) - F_0(x_1)), x_0 - x_1 \rangle &\geq& \langle (F_0(x_0) - F_0(x_1)), x_0 - x_1 \rangle \\ &\geq& \mu \| x_0 - x_1 \|^2
\end{eqnarray*} 
from which the theorem follows immediately $\qed$

Interestingly, in the above analysis, we did not assume any property of the intervened causal VI $F_1$, other than it has a solution (meaning that $F_1$ should be continuous). The following corollaries follow directly from Theorem~\ref{sens-thm1}. 
\begin{theorem}
\label{sens-thm2}
Given the original ``untreated" causal VI CVI$(F_0,K)$, where $F_0$ is strictly monotone, and the intervened ``treated" causal VI CVI$(F_1, K)$, where the intervened mapping $F_1$ is continuous, but not necessarily monotone, if $x_0$ and $x_1$ denote the solutions to the original ``untreated" and causally intervened CVI, where $x_0 \neq x_1$, then it follows that:
\begin{eqnarray}
\langle (F_1(x_1) - F_0(x_1), x_1 - x_0 \rangle < 0 \\
\langle (F_1(x_1) - F_0(x_0), x_1 - x_0 \rangle \leq 0
\end{eqnarray}
\end{theorem}
 Here, we are bounding the causal intervention effect of $F_1 - F_0$ of the ``treated" vs. ``untreated" operator, whereas previously in the convergence analysis of Equation~\ref{2stepalg}, we were trying to bound the same operator's effect on two different parameter values. To localize  causal interventions, we use the concept of a partitioned CVI. 
\begin{definition}
\label{partioned-cvi}
A partitioned CVI is defined as the causal  variational inequality problem of finding a vector $x^* = (x^*_1, \ldots, x^*_n) \in K \subset \mathbb{R}^n$ such that
\begin{equation*}
\langle E_w[F(x, \eta | \hat{w})], (x - y) \rangle \geq 0, \ \forall y \in K
\end{equation*}
where the function $F$ is partitionable function of order $m$, meaning that
\begin{equation*}
\langle E_w[F(x, \eta | \hat{w})], (x - y) \rangle = \sum_{i=1}^m \langle E_w[F_i(x, \eta | \hat{w})], (x_i - y_i) \rangle 
\end{equation*}
where each $F_i: K_i \subset \mathbb{R}^{n_i} \rightarrow \mathbb{R}^{n_i}$, with each $K_i$ being a convex domain such that $\prod_i K_i = K$.
\end{definition}
 The following theorem extends Theorem~\ref{sens-thm2} in showing that for partitionable CVI's, the effects induced by local causal interventions can be isolated.
\begin{theorem}
If a partitioned causal VI $CVI(F,K)$ is defined, where each component $F_i$ is a strongly monotone partitionable function,  $F^1_i$ denotes the causally intervened component function, and  $F^0_i$ is the ``untreated" function,  $x^0$ denotes the solution to the original ``untreated" CVI$(F_0,K)$ and $x^1$ denotes the solution to the ``treated" causally intervened CVI$(F_1,K)$ defined by the manipulated $F^1_i$ component functions, then 
\begin{equation}
\sum_{i=1}^m \langle (F^1_i(x^1_i) - F^0_i(x^1_i)), x^1_i - x^0_i \rangle  < 0
\end{equation}
\end{theorem}
{\bf Proof:} The proof follows readily from Theorem~\ref{sens-thm1}, Theorem~\ref{sens-thm2}, and Definition~\ref{partioned-cvi}. In particular, if $x^1$ is the solution of the intervened CVI$(F_1,K)$ and $x^0$ is the solution of the original CVI$(F_0,K)$, it follows that:
\[ \langle F_1(x^1) - F_0(x^1), x^1 - x^0 \rangle = \sum_{i=1}^m \langle F^1_i(x^1_i) - F^0_i(x^1_i), x^1_i - x^0_i \rangle \] 
Since the component functions $F_i$ are strongly monotone, the overall function $F$ is as well, and by applying Theorem~\ref{sens-thm1}, it follows that:
\[ \langle (F^1(x^1) - F^0(x^1)), x^1 - x^0 \rangle < 0 \] 
which immediately yields that
\[  \sum_{i=1}^m \langle (F^1_i(x^1_i) - F^0_i(x^1_i)), x^1_i - x^0_i \rangle < 0 \qed \]
If only a single component function $F_i$ is treated, then: 
\[  \langle (F^1_i(x^1_i) - F^0_i(x^1_i)), x^1_i - x^0_i \rangle < 0 \]
We can use these insights into designing an improved version of the two-step stochastic approximation algorithm given by Equation~\ref{2stepalg}. Instead of selecting random iterates to project on, we can instead prioritize those components $F^1_i$ that have been modified by the intervention.

%
%

\section{Summary and Related Work} 

Our paper proposed a novel variational framework for causal inference, by combining the ideas in causal inference \citep{rubin-book,pearl-book} and variational inequalities \citep{stampacchia,dafermos,nagurney:pdsbook,nagurney:vibook,facchinei-pang:vi}. We are not aware of previous work that explicitly connected these two longstanding research fields. 

While there has been significant work on analyzing causal effects in multiplayer games \citep{DBLP:journals/jmlr/BottouPCCCPRSS13,DBLP:conf/nips/ToulisP16,Wager2019ExperimentingIE}, our framework is significantly broader in scope and even in the special setting of games can be useful in generalizing previous work to nested levels of competing agents. Our causal modeling paradigm can be used in many other applications of VI's, including reinforcement learning. Our causal VI framework nicely complements previous work in RL using VI's \citep{bertsekas:vi-tr09,rl-lcp}, but significantly differs from previous work on combining reinforcement learning and structural causal models \citep{causal-rl}. Our causal VI formulation can also be profitably applied to VI formulations of deep learning models, such as generative adversarial networks \citep{gemp-gan}

\end{document}